\documentclass[12pt]{article}
\def\triangulo{\hbox{.\kern2pt.\kern-5pt\raise 4pt\hbox{.}}\kern 5pt}
\font\gordas = msbm10 at 12pt
\def\bbb#1{\hbox {{\gordas #1}}}
\def\a{{\bbb A}}
\def\erre{{\bbb R}}
\def\ce{{\bbb C}}
\def\o{{\bbb O}}

\def\ache{{\bbb H}}

\begin{document}
\begin{center}
{\large\bf Alternative elements in the
Cayley--Dickson algebras}\\[.5cm]
Guillermo Moreno\\
Departamento de Matem\'aticas\\
CINVESTAV-IPN, A.P. 14-7000\\
M\'exico, D.F. 07340
\end{center}

\begin{flushright}
To Jerzy Plebansky on his 75th. birthday
\end{flushright}
\vglue1cm
\noindent
{\bf Introduction.}
\vglue.5cm
$\a_n=\erre^{2^n}$ denotes the Cayley-Dickson algebra over $\erre$ the real numbers.

For $n=0,1,2,3\quad\a_n$ are the real, complex, quaternion and octonion numbers, denoted by $\erre, \ce, \ache$ and $\o$, respectively. 

These four algebras are normed,  
$||xy||=||x||||y||$  
and alternative,
$$x^2y=x(xy)\quad\mbox{\rm and}\quad xy^2=(xy)y$$
for all $x$ and $y$ in $\a_n$.

As is well known,  $\a_n$ is defined inductively by the Cayley-Dickson 
process:
\begin{eqnarray*}
x&=&(x_1,x_2)\quad y=(y_1,y_2)\quad\mbox{\rm in}\quad\a_n\times\a_n=\a_{n+1}\\
xy&=&(x_1y_1-\overline{y}_2x_2,y_2x_1+x_2\overline{y}_1)\quad\mbox{\rm and}\\
\overline{x}&=&(\overline{x}_1,-x_2).
\end{eqnarray*}
For $n\geq 4,\quad \a_n$ is a
neither normed nor alternative algebra, but $\a_n$ is {\it flexible}, i.e.
$$(xy)x=x(yx)\quad\mbox{\rm for all}\quad x \quad{\rm and}\quad y\quad{\rm in}
\quad \a_n$$
(see [3]).

In this paper we characterize the subset of $\a_n$ consisting of alternative elements, i.e. $\{a\in\a_n|a(ax)=a^2x$ for all $x$ in 
$\a_n\}$ and in terms  of  this characterization we 
``measure'' the failure of $\a_n\quad ( n\geq 4)$ of being normed .

Introducing the {\it associator} notation
$$(a,b,c):=(ab)c-a(bc).$$

We have that $a$ in $\a_n$ is alternative if and only if for all $x$ in $\a_n$
$$(a,a,x)=0$$

$\a_n$ is flexible if and only if $(a,x,a)=0$ for all $a$ and $x$ in $\a_n$.

Now $\{e_0,e_1,e_2,\ldots,e_{2^n-1}\}$ denote the canonical basis in $\a_n$.

$\a_n={\rm Span}\{e_0\}\oplus{\rm Span}\{e_1,e_2,\ldots, e_{2^n-1}\}:=\erre e_0\oplus Im(\a_n)$ is the 
canonical splitting into real and pure imaginarie part in $\a_n$.

In $\S$2 we prove that $(a,x,b)=0$ for all $x$ in $\a_n$ if and only if $a$ and $b$ have linearly dependent pure imaginarie parts.This statement is known as the Yui Conjecture.

We present a different proof of this fact,from [1] and [2].

Using this, we characterize the alternative elements in $\a_n$ for all $n\geq 4$.

Here is an outline of the main ideas.

First of all we notice that if $L_a:\a_n\rightarrow\a_n$ denotes the left multiplication by $a$ in $\a_n$ then $L_a$ is a linear transformation and that for $a$ pure element,  
$L_a$ is skew--symmetric so $L^2_a$ is symmetric non-positive definite [3].

But $a$ is alternative if and only if its imaginary part is alternative and for $a$ pure alternative
 element $L^2_a=a^2I=-||a||^2I$ i.e.;  $L^2_a$ has all its eigenvalues equal to $-||a||^2$. 
By direct calculation we will see that for $a$ and $b$ pure elements in $\a_n$%
$$L^2_{(a,b)}:\a_{n+1}\rightarrow\a_{n+1}$$
is given by
$$L^2_{(a,b)}(x,y)=({\cal A}(x)-S(y),{\cal A}(y)+S(x))$$
where 
$${\cal A}=L^2_a+R^2_b:\a_n\rightarrow\a_n$$
and
$$S=(a,-,b)=R_bL_a-L_aR_b:\a_n\rightarrow\a_n$$
where $R_b$ is the right multiplication by $b$ (see Lemma 3.2) So if $(a,b)$ is alternative in $\a_{n+1}$ 
then $({\cal A}(x)-S(y),{\cal A}(y)+S(x))=(a,b)^2(x,y)$ for all $x$ and $y$
then $S(x)=0$ for all $x$ in $\a_n$ 
then $a$ and $b$ are linearly dependent and ${\cal A}(x)=(a^2+b^2)x$ 
for all $x$ so $a$ and $b$ are alternatives (see Theorem 3.3)

Thus alternative elements in $\a_{n+1}$ are ``made from'' the alternative 
elements in $\a_n$. Despite  of the fact that in $\a_3$ all element is 
alternative, in $\a_4$ they form a ``very small'' subset, this set remains 
``constant'' during the doubling process.

In $\S 4$ we define 
$a$ in $\a_n$ is {\it strongly alternative} if $(a,a,x)=(a,x,x)=0$ for all $x$ 
in $\a_n$. Then we prove the set of strongly alternative elements, form even a smaller subset, 
namely 
$\{re_0+s\widetilde{e}_0|r$ and $s$ in $\erre\}$ where 
$\widetilde{e}_0=e_{2^{n-1}}$.

In $\S 5$ we study the properties of being alternative and strongly alternative locally, i.e.; between two elements in $\a_n$.

This give us a criterion to know which elements generate associative and alternative subalgebras inside of $\a_n\quad n\geq 4$.

In this paper, sequel of [3], we use that the canonical basis $\{e_0,e_1\ldots, e_{2^n-1}\}$ consists of alternative elements and that the Euclidean structure in $\erre^{2^n}$ and the $C-D$ algebra structure are related by 
$$2\langle x,y\rangle=x\overline{y}+y\overline{x}$$%
$||x||^2=x\overline{x}$ for all $x$ and $y$ in $\a_n$. (See [4]).
\vglue1cm
\noindent
{\bf \S 1.  Pure and doubly pure elements.}
\vglue.5cm

For $a=(a_1,a_2)\in\a_{n-1}\times \a_{n-1}=\a_n$ {\it the trace} is
\begin{eqnarray*}
\lefteqn{t_n:\a_n\rightarrow\a_0=\erre}\\
t_n(a)&=&a+\overline{a}.
\end{eqnarray*}
{\bf Definition.} We say that $a\in\a_n$ is {\it pure} if
$$t_n(a)=0\quad i.e.\quad \overline{a}=-a.$$
Notice that $a$ is pure in $\a_n$ if and only if $a_1$ is pure in $\a_{n-1}$.

Let $\{e_0,e_1,\ldots,e_{2^n-1}\}$ be the canonical basis in $\erre^{2^n}$ 
so $e_0=(1,0,\ldots,0)$ is the unit in the algebra $\a_n$.
\vglue.5cm
\noindent
{\bf Notation.} $_0\a_n:=\{a\in\a_n|a$ is pure $\}=\{e_0\}^\perp$.

For $a=(a_1,a_2)\in\a_{n-1}\times\a_{n-1}=\a_n$ we denote
$$\widetilde{a}=(-a_2,a_1)\in\a_{n-1}\times\a_{n-1}=\a_n.$$
Note that $||a||^2=||\widetilde{a}||^2=||a_1||^2+||a_2||^2$. 

In terms of the decomposition $\a_n=\a_{n-1}\times\a_{n-1}$%
\begin{eqnarray*}
e_0=(e_0,0),e_1&=&(e_1,0),\ldots, e_{2^{n-1}-1}=(e_{2^{n-1}-1},0),\\
e_{2^{n-1}}=(0,e_0),e_{2^{n-1}+1}&=&(0,e_1),\ldots, 
e_{2^n-1}=(0,e_{2^{n-1}-1}),
\end{eqnarray*}
we have that $\widetilde{e}_i=e_{2^{n-1}+i}$ for $0\leq i\leq 2^{n-1}$ and
$\widetilde{e}_i=-e_{i-2^{n-1}}$ for 
$$2^{n-1}< i\leq 2^{n}-1,$$ 
so
$\widetilde{e}_0=e_{2^{n-1}}$ and $\widetilde{a}=a\widetilde{e}_0$.
\vglue.5cm
\noindent
{\bf Definition.} An element $a=(a_1,a_2)\in\a_{n-1}\times \a_{n-1}=\a_n$ {\it is doubly pure} if both coordinates $a_1$ and $a_2$ are pure elements in $\a_{n-1}$ i.e. $t_{n-1}(a_1)=t_{n-1}(a_2)=0$. Notice that $a$ in $\a_n$ is doubly pure if and only if $a$ and $\widetilde{a}$ are pure elements in $\a_n$.

Notice that for any $a$ pure element in $\a_n$,  $a=rc+s\widetilde{e}_0$ for $r$ and $s$ real numbers and $c$ doubly pure element in $\a_n$.
\vglue.5cm
\noindent
{\bf Notation.} $\widetilde{\a}_n:=\{a\in\a_n| a$ is doubly pure $\}$. 

Notice that 
$\widetilde{\a}_n=\{e_0,\widetilde{e}_0\}^\perp =
_o\!\!\a_{n-1}\times _o\!\!\a_{n-1}\subset\a_n$.

Now the Euclidean product in $\a_n=\erre^{2^n}$ is given by 
$$2\langle x,y\rangle=t_n(x\overline{y})=x\overline{y}+y\overline{x}$$ 
and the Euclidean norm by $||x||^2=x\overline{x}$. For $a=(a_1,a_2)\in\a_{n-1}\times\a_{n-1}=\a_n$ doubly pure we have that $a\in\{e_0,\widetilde{e}_0\}^\perp$ so $\widetilde{a}=a\widetilde{e}_0=-\widetilde{e}_0a$ and  $\widetilde{a}a+a\widetilde{a}=(a\widetilde{e}_0)a-a(\widetilde{e}_0a)=(a,\widetilde{e}_0,a)=0$ (by flexibility) then $\widetilde{a}\perp a$.
\vglue.5cm
\noindent
{\bf Lemma 1.1.} Let $a=(a_1,a_2)$ and $x=(x_1,x_2)$ be elements in $\a_{n-1}\times\a_{n-1}=\a_n$.

If $x$ is doubly pure then $\widetilde{a}x=-\widetilde{ax}$.
\vglue.5cm
\noindent
{\bf Proof.} Define $c=a_1x_1+x_2a_2$ and $d=x_2a_1-a_2x_1$ in $\a_{n-1}$ 
so
$ax=(a_1x_1+x_2a_2,x_2a_1-a_2x_1)=(c,d)$ and $\widetilde{a}
x=(-a_2x_1+x_2a_1,-x_2a_2-a_1x_1)=(+d,-c)$ and $\widetilde{ax}=(-d,c)$ 
so $\widetilde{a}x=-\widetilde{ax}$. \hfill Q.E.D.
\vglue.5cm
\noindent
{\bf Corollary 1.2.} For $a$ and $x$
doubly pure elements 
in $\a_n$
we have that

\begin{enumerate}
\item[1)] $\widetilde{a}x+\widetilde{x}a=0$ if and only if $a\perp x$.
\item[2)] $ax-\widetilde{x}\widetilde{a}=0$ if and only if $\widetilde{a}\perp x$.
\item[3)] $\widetilde{a}x=0$ if and only if $ax=0$.
\end{enumerate}
{\bf Proof.}

\begin{enumerate}
\item[1)] $a\perp x$ if and only if $ax=-xa$ $\Leftrightarrow\widetilde{ax}=-\widetilde{xa}\Leftrightarrow\widetilde{a}x=-\widetilde{x}a$.
\item[2)]
$\widetilde{a}\perp x\Leftrightarrow\widetilde{a}x=-x\widetilde{a}\Leftrightarrow\widetilde{\widetilde{a}x}=-\widetilde{x\widetilde{a}}\Leftrightarrow ax=\widetilde{x}\widetilde{a}$.
\item[3)] $0=ax\Leftrightarrow 0=\widetilde{ax}\Leftrightarrow 0=\widetilde{a}x$. \hfill Q.E.D.
\end{enumerate}
{\bf Corollary 1.3.} For $0\neq a$  in $\a_n$ doubly pure element and 
 $n\geq 3$ the vector subspace of $\a_n$ generated by $\{e_0,\widetilde{a},a,\widetilde{e}_0\}$ is a copy of the quaternions $\a_2$.
\vglue.5cm
\noindent
{\bf Proof.} Since 
$a\in\{e_0,\widetilde{e}_0\}^\perp$ then 
$\widetilde{a}\in\{e_0,\widetilde{e}_0\}^\perp$ and $a\perp\widetilde{a}$ 
so $\{e_0,\widetilde{a},a,\widetilde{e}_0\}$ is an orthogonal set of four vectors in 
$\a_n$ for $n\geq 3$. (we denote it by $\ache_a$). Now we suppose that 
$||a||=1$ otherwise we take ${a\over||a||}$. From Lemma 1.1 and Corollary 1.2 
we have the following quaternion multiplication table.
$$
\begin{tabular}{l|llll}
&$e_0$ &$\widetilde{a} $&$a$&$\widetilde{e}_0$\\ \hline
$e_0$&$e_0$&$\widetilde{a}$ &$a$&$\widetilde{e}_0$\\
$\widetilde{a}$&$\widetilde{a}$&$-e_0$&$-\widetilde{e}_0$&$-a$\\
$a$&$a$&$\widetilde{e}_0$&$-e_0$&$\widetilde{a}$\\
$\widetilde{e}_0$&$\widetilde{e}_0$&$a$&$-\widetilde{a}$&$-e_0$
\end{tabular}
$$
here $e_0\leftrightarrow 1; \widehat{i}\leftrightarrow\widetilde{a};\hat{j}
\leftrightarrow a$ and $\widehat{k}\leftrightarrow\widetilde{e}_0$. 
\hfill Q.E.D.
\vglue.5cm
\noindent
{\bf Proposition 1.4.} For $a$ and $b$ doubly pure elements in $\a_n$ for $n\geq 3$ we have:

\begin{enumerate}
\item[1)] $\widetilde{a}b=a\widetilde{b}$ if and only if $a\perp b$ and $\widetilde{a}\perp b$.
\item[2)] If $(a,\widetilde{e}_0, b)=0$ then $ab=re_0+s\widetilde{e}_0$ for $r$ and $s$ in $\erre$.
\end{enumerate}
{\bf Proof.} 
Let $a=(a_1,a_2)$ and $b=(b_1,b_2)$ be in $_o\a_{n-1}\times _o\a_{n-1}$. If 
$r=a_1b_1+b_2a_2$ and $s=b_2a_1-a_2b_1$ in 
$\a_{n-1}$ then $ab=(r,s)\in\a_{n-1}\times\a_{n-1}=\a_n$. 
So $\widetilde{ab}=(-s,r)$ and $\widetilde{a}b=(s,-r)$. On the other hand 
$a\widetilde{b}=(a_1,a_2)(-b_2,b_1)=(-a_1b_2+b_1a_2,b_1a_1+a_2b_2)$ so $a\widetilde{b}=(-\overline{s},\overline{r})$. Therefore $\widetilde{a}b=a\widetilde{b}\Leftrightarrow s=-\overline{s}$ and $\overline{r}=-r\Leftrightarrow t_{n-1}(s)=0$ and $t_{n-1}(r)=0\Leftrightarrow t_n(ab)=0$ and $t_n(\widetilde{a}b)=0\Leftrightarrow a\perp b$ and $\widetilde{a}\perp b$ so we prove 1).

To prove 2) we see that $(a,\widetilde{e}_0,b)=\widetilde{a}b+a\widetilde{b}$ so if $0=(a,\widetilde{e}_0,b)$ in $\a_n$ then
$$(0,0)=(s,-r)+(-\overline{s},\overline{r})=(s-\overline{s},-r+\overline{r})\quad\mbox{\rm in}\quad \a_{n-1}\times\a_{n-1}$$
so $\overline{s}=s$ and $\overline{r}=r$ and $s$ and $r$ are real 
numbers therefore $ab=re_0+s\widetilde{e}_0$. 

\hfill Q.E.D.
\vglue.5cm
\noindent
{\bf Corollary 1.5.} For $a$ and $b$ doubly pure elements in $\a_n$ for $n\geq 3$, we have that
$$
-(\widetilde{e}_0,a,b)=(a,\widetilde{e}_0,b)=\left\{
\begin{array}{lll}
0&if&b\in\ache_a\\
2\widetilde{a}b&if&b\in\ache^\perp_a
\end{array}
\right.
$$
{\bf Proof.} Since $\ache_a$ is associative $-(\widetilde{e}_0,a,b)=(a,\widetilde{e}_0,b)=0$ for $b\in\ache_a$. If $b\in\ache^\perp_a$ then $a\perp b$ and $\widetilde{a}\perp b$ and by Proposition 1.4 (2)
$$(a,\widetilde{e}_0,b)=\widetilde{a}b+a\widetilde{b}=2\widetilde{a}b=-
(\widetilde{e}_0a)b+\widetilde{e}_0(ab)=-(\widetilde{e}_0,a,b).$$ 

\hfill Q.E.D.
\vglue.5cm
\noindent
{\bf Remark.} Notice that also we prove that for $a$ and $b$ doubly pures, if $b\in\ache_a$ then $ab=re_0+s\widetilde{e}_0$ for $r$ and $s$ in $\erre$.
\vglue.5cm
\noindent
{\bf Lemma 1.6.} Let $a$ be a doubly pure element in $\a_n$ for $n\geq 3$ we have that
$$
(\widetilde{a},x,a)=\left\{
\begin{array}{lll}
0&if&x\in\ache_a\\
-2a(a\widetilde{x})&if&x\in\ache^\perp_a
\end{array}
\right.
$$
{\bf Proof.} Since $\ache_a=\langle\{e_0,\widetilde{a},a,\widetilde{e}_0\}\rangle$ is associative
$$(\widetilde{a},x,a)=0\quad\mbox{\rm for }\quad x\in\ache_a.$$
Suppose that $0\neq x\in\ache^\perp_a$ then $x$ is doubly pure and $a\perp x$ and $\widetilde{a}\perp x$. By Proposition 1.4 (1) and flexibility $(\widetilde{a}x)a=(a\widetilde{x})a=a(\widetilde{x}a)=-a(a\widetilde{x})$ (recall that $a\perp\widetilde{x}$). Since $a\perp x$ then $ax$ is pure and $\widetilde{a}\perp x$ implies that $ax$ is doubly pure then applying Proposition 1.4 (1) to $a$ and $xa$ we have that
$$\widetilde{a}(xa)=a(\widetilde{xa})=-a(\widetilde{x}a)=a(a\widetilde{x})$$
because $a\perp xa$ and $\widetilde{a}\perp xa$ (recall that right multiplication by any pure element is a skew symmetric linear transformation). Therefore $$(\widetilde{a},x,a)=(\widetilde{a}x)a-\widetilde{a}(xa)=-2a(a\widetilde{x})$$. \hfill Q.E.D.
\vglue1cm
\noindent
{\bf \S 2. Proof of the Yui's Conjecture.}
\vglue.5cm
In this section we give an affirmative answer to the following question: 

By flexibility we have that: If $a$ and $b$ are linearly dependent in $_o\!\a_n$ with $n\geq 3$ then $(a,x,b)=0$ for all $x$ in $\a_n$.

{\bf Question:} Is the converse true?.

We first show the case for $a$ and $b$ doubly pure elements and then we proceed with the general case:
\vglue.5cm
\noindent
{\bf Lemma 2.1.} Let $a\in\a_n$ be a doubly pure element with $n\geq 4$. 

If $ax=0$ for all $x\in\ache^\perp_a\quad$ then $a=0$.
\vglue.5cm
\noindent
{\bf Proof.} Let `$\varepsilon$' denote the basic element $e_{2^{n-2}}$ in $\a_n$ i.e. $\varepsilon=e_4$ in $\a_4\quad\varepsilon=e_8$ in $\a_5,\ldots,$ etc. Therefore $\varepsilon$ is an alternative element so $\varepsilon$ can't be a zero divisor. Suppose that $ax=0$ for all $x\in\ache^\perp_a\subset\a_n$ and $n\geq 4$. If $\varepsilon\in\ache^\perp_a$ then $a\varepsilon=0$ and $a=0$. If $\varepsilon\in\ache_a$ then $\ache_a=\ache_\varepsilon$ and $ax=0$ if and only if $\varepsilon x=0$ i.e. $x=0$ and $\ache^\perp_a=\{0\}$ and $\a_n=\ache_a$ but $n\geq 4$ and $2^n>4=dim\ache_a$ (contradiction). Let's suppose that $\varepsilon=(\varepsilon '+\varepsilon '')\in\ache_a\oplus\ache^\perp_a$ 
with $\varepsilon ''\neq 0$ so $a\varepsilon ''=0$.

Therefore $a\varepsilon=a\varepsilon '+a\varepsilon ''=a\varepsilon 
'+0=a\varepsilon '\in\ache_a$. 

Since $a$ is doubly pure, $a\varepsilon\perp a$ and $a\varepsilon\perp\widetilde{a}$ 
so
we have  $a\varepsilon\in\ache_a$ implies that $a\varepsilon=re_0+s\widetilde{e}_0$ for some $r$ and $s$ in $\erre$, and $a(a\varepsilon)=ra+s\widetilde{a}$.

 On the other hand $a(a\varepsilon)=a(a\varepsilon ')=a^2\varepsilon '$ because $\ache_a$ is associative so $-||a||^2\varepsilon '=a^2\varepsilon=ra+s\widetilde{a}$ and $-||a||^2\varepsilon '\varepsilon ''=(ra+s\widetilde{a})\varepsilon ''=ra\varepsilon ''+s\widetilde{a}\varepsilon ''=0+0=0$ (Lemma 1.1) and $\varepsilon '\varepsilon ''=0$ unless $a=0$.If $a\neq 0$ then $\varepsilon '\varepsilon=\varepsilon '(\varepsilon '+\varepsilon '')=\varepsilon'^2+0=-||\varepsilon '||^2\in\erre$ so $\varepsilon$ and $\varepsilon '$ are linearly dependent and $\varepsilon\in\ache_a$. But this is impossible,so $a=0$. \hfill Q.E.D.

Notice that in this proof we use only the fact that $\varepsilon$ is an alternative element of norm one.

\vglue.5cm
\noindent
{\bf Theorem 2.2.} For $a$ and $b$ non-zero doubly pure elements in $\a_n$ and $n\geq 4$ we have: If $(a,y,b)=0$ for all $y\in\a_n$ then $a$ and $b$ are linearly dependent.
\vglue.5cm
\noindent
{\bf Proof.} We proceed by contradiction. Suppose that $a$ and $b$ are linearly independent. Without loss of generality we can suppose that $a$ is orthogonal to $b$ because, by flexibility, $(a,y,b)=(a,y,b-ra)$ for all $y$ and $r={\langle b,a\rangle\over\langle a,a\rangle}$. On the other hand if $(a,\widetilde{e}_0,b)=0$ then by (Proposition 1.4 (2)) $ab=pe_0+q\widetilde{e}_0$ for $p$ and $q$ in $\erre$, but $a\perp b$ so $p=0$ and $ab=q\widetilde{e}_0$ then $\widetilde{a}$ and $b$ are linearly dependent. Now 
$(a,x,\widetilde{a})=-(\widetilde{a},x,a)=0$ for all $x\in\a_n$ implies that $a(a\widetilde{x})=0$ for all $x\in\ache^\perp_a$ (Lemma 1.6) so $ax=0$ for all $x\in\ache^\perp_a$ (recall that $KerL_a = KerL^2_a$). By Lemma 2.1 $a=0$ which is a contradiction. \hfill Q.E.D.
\vglue.5cm
Now we proceed with the general case: Since any associator with one entrie equal to $e_0$ automatically vanishes we have to prove that if $(\alpha,x,\beta)=0$ for all $x\in\a_n$ then $\alpha$ and $\beta$ are linearly dependent for $\alpha$ and $\beta$ pure. 

But  $\alpha=a+p\widetilde{e}_0$ and $\beta=b+q\widetilde{e}_0$ where $a$ and $b$ are doubly pure elements and $p$ and $q$ are real numbers.

From now on we suppose that:
  $(\alpha,x,\beta)=0\quad\forall x\in\a_n\quad(n\geq 4)$  

Suppose that $b=0$.\\
 Thus $(a+p\widetilde{e}_0,x,\widetilde{e}_0)=(a,x,\widetilde{e}_0)+p(\widetilde{e}_0,x,\widetilde{e}_0)=(a,x,\widetilde{e}_0)$ but $(a,x,\widetilde{e}_0)=\widetilde{ax}-a\widetilde{x}=
-\widetilde{a}x-a\widetilde{x}=-2\widetilde{a}x$ for all $x\in\ache^\perp_a$.

So $(\alpha,x,\beta)=0$ if and only if $ax=0$ for all $x\in\ache^\perp_a$ and by Lemma 2.1  $a=0$ and  $\alpha$ and $\beta$ are linearly dependent. 
The argument is similar for  $a=0$.
 
Suppose that $a\neq 0$ and $b\neq 0$%
$$\qquad
(\alpha,x,\beta)=(a+p\widetilde{e}_0,x,b+q\widetilde{e}_0)=(a,x,b)+(qa-pb,x,\widetilde{e}_0). \qquad(*)$$
Put  $x=\widetilde{e}_0$
$$(\alpha,\widetilde{e}_0,\beta)=(a,\widetilde{e}_0,b)+(qa-pb,\widetilde{e}_0,\widetilde{e}_0)=(a,\widetilde{e}_0,b)+0.$$
So,  $(\alpha,\widetilde{e}_0,\beta)=(a,\widetilde{e}_0,b)=0$ then by 
Proposition 1.4 (2) 
$$ab=re_0+s\widetilde{e}_0$$   for $r$  and $s$ real numbers.
 Put  $x=a$%
$$(\alpha,a,\beta)=(a,a,b)+(qa-pb,a,\widetilde{e}_0)=(a,a,b)+q(a,a,\widetilde{e}_0)+p(b,a,\widetilde{e}_0)$$
but $(a,a,\widetilde{e}_0)=0$ because $\ache_a$ is associative and 
$$(b,a,\widetilde{e}_0)=\widetilde{ba}-b\widetilde{a}=-\widetilde{b}a-b\widetilde{a}=-(b,\widetilde{e}_0,a)=(a,\widetilde{e}_0,b)=0$$
Therefore  $(\alpha,a,\beta)=(a,a,b)=0$ and
$$a^2b=a(ab)=a(re_0+s\widetilde{e}_0)=ra+s\widetilde{a}$$.
Since $a\neq 0\quad b=ua+v\widetilde{a}$ where $u={r\over a^2}, v={s\over a^2}$and substituting in (*) above we have
$$
0=(\alpha,x,\beta)=(a,x,b)+(qa-pb,x,\widetilde{e}_0)
$$
\begin{eqnarray*}
&=&(a,x,ua+v\widetilde{a})+(qa-p(ua+v\widetilde{a}),x,\widetilde{e}_0)\\
&=&u(a,x,a)+v(a,x,\widetilde{a})+(q-pu)(a,x,\widetilde{e}_0)-pv(\widetilde{a},x,\widetilde{e}_0)\\
&=&0+v(a,x,\widetilde{a})+(q-pu)(a,x,\widetilde{e}_0)+pv(\widetilde{a},x,
\widetilde{e}_0).
\end{eqnarray*}
Now $(a,x,\widetilde{a})=(a,x,\widetilde{e}_0)=(\widetilde{a},x,\widetilde{e}_0)=0$ for $x\in\ache_a$ \\
and
$(a,x,\widetilde{a})=-(\widetilde{a},x,a)=2a(a\widetilde{x})$ for 
$x\in\ache^\perp_a$ (Lemma 1.6) 
\begin{eqnarray*}
(a,x,\widetilde{e}_0)&=&\widetilde{ax}-a\widetilde{x}=-\widetilde{a}x-a\widetilde{x}=-2\widetilde{a}x\quad{\rm for}\quad x\in\ache^\perp_a\\
(\widetilde{a},x,\widetilde{e}_0)&=&(\widetilde{\widetilde{a}x})-
\widetilde{a}\widetilde{x}=ax-xa=2ax\quad{\rm for}\quad x\in\ache^\perp_a.
\end{eqnarray*}
Now multiplication by $\widetilde{e}_0$ and `$a$' are skew-symmetric linear transformations so $\{-\widetilde{a}x=\widetilde{ax},ax,a(a\widetilde{x})\}$ is an orthogonal subset.

By Lemma 2.1 $ax=0, \widetilde{a}x=0$ and $a(a\widetilde{x})=0$ for all $x\in\ache^\perp_a$ only if $a=0$.Therefore $0=(\alpha,x,\beta)$ for all $x$ in $\a_n$ implies that:
$$v=0\quad q-pu=0\quad{\rm and}\quad pv=0$$ 
then $s=0\quad q=pu$. 

To finish we argue as follows:

Since $b=ua$ then $b+q\widetilde{e}_0=ua+pu\widetilde{e}_0=u(a+p\widetilde{e}_0)$, i.e. $\beta=u\alpha$ and $\alpha$ and $\beta$ are linearly dependent, 
so we proved                 :
\vglue.5cm
\noindent
{\bf Theorem 2.3.} If $a$ and $b$ are non-zero elements in $\a_n$ for $n\geq 4$ such that $(a,x,b)=0$ for all $x$ in $\a_n$ then $a$ and $b$ 
have linearly dependent pure parts.
\vglue1cm
\noindent
{\bf \S 3. Alternative elements in $\a_n\quad n\geq 4$.}
\vglue.5cm
\noindent
{\bf Definition.} $a\in\a_n$ is an 
{\it alternative element} if $(a,a,x)=0$ for all $x$ in $\a_n$.

It is known [4] that the elements in the canonical basis are alternative.

Clearly a scalar multiple of an alternative element is alternative but the sum of two alternative elements is not necessarily alternative. 

Because the associator symbol vanish if one of the entries is real, i.e. belongs to $\erre e_0$, then an element is alternative if and only if its pure (imaginary) part is alternative,
Therefore  we need to characterize the pure (non-zero) alternative elements.
\vglue.5cm
\noindent
{\bf Proposition 3.1.} Let $a\in\a_n\quad (n\geq 4)$ be a pure element with 
$$a=rc+s\widetilde{e}_0$$%
$r\neq 0$ and $s$ in $\erre$ and $c$ doubly pure. Then $a$ is alternative if and only if $c$ is alternative.

(We define `$c$' as the {\it doubly pure part } of $a$).
\vglue.5cm
\noindent
{\bf Proof.} For all $x$ in $\a_n$
\begin{eqnarray*}
(a,a,x)&=&(rc+s\widetilde{e}_0,rc+s\widetilde{e}_0,x)\\
&=&r^2(c,c,x)+s^2(\widetilde{e}_0,\widetilde{e}_0,x)+rs(c,\widetilde{e}_0,x)+rs(\widetilde{e}_0,c,x).
\end{eqnarray*}
But $\widetilde{e}_0$ is alternative so $(\widetilde{e}_0,\widetilde{e}_0,x)=0$ and by Corollary 1.5 $(c,\widetilde{e}_0,x)+(\widetilde{e}_0,c,x)=0$ so $(a,a,x)=r^2(c,c,x)$. 

If $r=0$ then $a=s\widetilde{e}_0$ and $a$ is alternative.

If $r\neq 0$ then $(a,a,x)=0$ if and only if $(c,c,x)=0$ for all $x\in\a_n$. \hfill Q.E.D.
\vglue.5cm
Because of Proposition 3.1 we will focus on the doubly pure alternative elements in $\a_n$ for $n\geq 4$.

For $a$ and $b$ pure elements in $\a_n$ and  $n\geq 3$ we denote $L_a,R_b:\a_n\rightarrow\a_n$ the left and right multiplication by $a$ and $b$ and $L_{(a,b)}:\a_{n+1}\rightarrow\a_{n+1}$ the left multiplication by $(a,b)$ in $\a_{n+1}$. Note that by flexibility $(a,a,x)=-(x,a,a)$ and $a(ax)=(xa)a$ 
i.e. $L^2_a=R^2_a$.
\vglue.5cm
\noindent
{\bf Notation.} For $a$ and $b$ fixed pure elements
$${\cal A}:=L^2_a+R^2_b$$%
$$S:=(a,-,b)=R_bL_a-L_aR_b=[R_b,L_a].$$
{\bf Lemma 3.2.} For $a$ and $b$ pure elements in $\a_n\quad n\geq 3$.
$$L^2_{(a,b)}(x,y)=({\cal A}(x)-S(y),{\cal A}(y)+S(x))$$
for $(x,y)\in\a_n\times\a_n=\a_{n+1}$.
\vglue.5cm
\noindent
{\bf Proof.} (Direct calculation).

$(a,b)(x,y)=(ax-\overline{y}b,ya+b\overline{x})$
and
\begin{eqnarray*}
(a,b)[(a,b)(x,y)]&=&(a,b)(ax-\overline{y}b,ya+b\overline{x})\\
&=&(a(ax-\overline{y}b)-(\overline{ya+b\overline{x}})b,(ya+b\overline{x})a+b(\overline{ax-\overline{y}b}))\\
&=&(a(ax)-a(\overline{y}b)-(\overline{ya})b-(x\overline{b})b,(ya)a+(b\overline{x})a+b(\overline{ax})+b(by))\\
&=&(a(ax)+(xb)b+(a\overline{y})b-a(\overline{y}b),a(ay)+(yb)b+(b\overline{x})a-b(\overline{x}a)
)\\
&=&(L^2_a(x)+R^2_b(x)+S(\overline{y}),L^2_a(y)+R^2_b(y)-S(\overline{x}))\\
&=&({\cal A}(x)-S(y),{\cal A}(y)+S(x)).
\end{eqnarray*}
because $(a,\overline{x},b)=-(a,x,b)$ for all $x$, $\overline{a}=-a$,$\overline{b}=-b$ and $b(by)=(yb)b$.

\hfill Q.E.D.
\vglue.5cm
\noindent
{\bf Remark.} Notice if we 
 interchange the role of $a$ and  $b$ we have that
$$L^2_{(b,a)}(x,y)=({\cal A}(x)+S(y), {\cal A}(y)-S(x))$$
by flexibility.
\vglue.5cm
\noindent
{\bf Theorem 3.3.} For $a$ and $b$ pure elements in $\a_n$, for 
$n\geq 3$, consider the following statements:

\begin{itemize}
\item[(i)] $(a,b)$ is an alternative element in $\a_{n+1}$.
\item[(ii)] $a$ and $b$ are alternative elements in $\a_n$.
\item[(iii)] $a$ and $b$ are linearly dependent in $\a_n$.
\end{itemize}

Then (i) if and only if (ii) and (iii).
\vglue.5cm
\noindent
{\bf Proof.} Suppose that $(a,b)$ is alternative in $\a_{n+1}$.
\vglue.5cm
Then
$L^2_{(a,b)}(x,y)=({\cal A}(x)-S(y),{\cal A}(y)+S(x))=(a,b)^2(x,y)$ 
by lemma 3.2.

So
\begin{eqnarray*}
{\cal A}(x)-S(y)&=&(a^2+b^2)x\\
{\cal A}(y)+S(x)&=&(a^2+b^2)y
\end{eqnarray*}
for all $x$ and $y$ in $\a_n$.

Let's put $x=y$ and substract,so 
$$S(x)=0\quad {\rm for\,\,\, all}\quad x\in\a_n.$$
By  Theorem 2.3 we have that 
$a$ and $b$ are linearly dependent. Now 
${\cal A}(x):=(L^2_a+R^2_b)(x)=(a^2+b^2)x$ and $a$ and $b$ 
linearly dependent implies that 
$a(ax)=a^2x$ and $(xb)b=b^2x$, i.e. $a$ and $b$ are alternative in 
$\a_n$. Conversely if $a$ and $b$ are linearly dependent, then 
$S(x):=(a,x,b)=0$ for all $x\in\a_n$ and $L^2_a$ and $R^2_b$
are multiples  of each  other. So 
$L^2_{(a,b)}(x,y)=({\cal A}(x),{\cal A}(y))=((L^2_a+R^2_b)(x),(L^2_a+R^2_b)(y))$  but $a$ and $b$ are also alternative so $L^2_a(x)=a^2x$ and $R^2_b(x)=b^2x$ for all $x$ in $\a_n$. So $L^2_{(a,b)}(x,y)=(a^2+b^2)(x,y)$ for all $(x,y)$ 
in $\a_{n+1}$ and $(a,b)$ is alternative in $\a_{n+1}$.

\hfill Q.E.D.
\vglue.5cm
\noindent
{\bf Notation.} 
\begin{eqnarray*}
Alt_n&=&\{a\in\a_n|a\quad{\rm is\,\,\, alternative}\}\\
Alt_n^0&=&\{a\in _o\!\!\a_n|a\quad{\rm is\,\,\, alternative}\}\\
\widetilde{Alt_n}&=&\{a\in\widetilde{\a}_n|a\quad{\rm is\,\,\, alternative}\}
\end{eqnarray*}
i.e. while $Alt_n$ denote the subset of $\a_n$ consisting of alternative elements, $Alt^0_n$ and $\widetilde{Alt}_n$ denote the pure and doubly pure alternative elements respectively in $\a_n$.

So by the decompositions
$$\a_n=_o\!\!\a_n\oplus\erre e_0=(\widetilde{\a}_n\times\erre\widetilde{e}_0)\times\erre e_0$$
we have 
$Alt_n\cong Alt^0_n\times\erre$ and $\widetilde{Alt}_n\times\erre=Alt^0_n$. Now $Alt_3=\a_3, Alt^0_3= _o\a_3$ and $\widetilde{Alt}_3=\widetilde{\a}_3$ so $\widetilde{Alt}_4=\{(ra,sa)\in\a_3\times\a_3|r$ and $s$ in $\erre\}$ and 
$\widetilde{Alt}_{n+1}=\{(ra,sa)\in Alt^0_n\times Alt^0_n|r\quad
{\rm  and}\quad s\quad{\rm  in}\quad \erre\}$ 
i.e.
$\widetilde{Alt}_{n+1}$ 
is a 'cone' on $Alt^0_n$ with an extra point $(r,s)=(0,0)$.
\vglue1cm
\noindent
{\bf \S 4. Strong alternativity}
\vglue.5cm
\noindent
{\bf Definition.} $a$ in $\a_n$ with $n\geq 3$ is{\ strongly alternative} if it is 
an
alternative element (i.e. $(a,a,x)=0$ for all $x$ in $\a_n$) 
and also $(a,x,x)=0$ for all $x$ in $\a_n$.
\vglue.5cm
\noindent
{\bf Example.} $a=e_1$ is alternative element but $a=e_1$ is not strongly alternative element in $\a_4$. Take $x=e_4+e_{15}$ in $\a_4$ so
\begin{eqnarray*}
(x,x,a)&=&(e_4+e_{15},e_4+e_{15},e_1)\\
&=&(e_4,e_4,e_1)+(e_{15},e_{15},e_1)+(e_4,e_{15},e_1)+(e_{15},e_4,e_1)\\
&=&0+0+(e_4e_{15}+e_{15}e_4)e_1-(e_4(e_{15}e_1)+e_{15}(e_4e_1))\\
&=&0+0+0-e_4e_{14}+e_{15}e_5\\
&=&e_{10}+e_{10}\\
&=&2e_{10}.
\end{eqnarray*}
Therefore among the alternative elements there are ``few'' strongly alternative elements.
\vglue.5cm
\noindent
{\bf Example.} $\widetilde{e}_0$ in $\a_n$ is strongly alternative for $n\geq 3$.

If $x$ is doubly pure element in $\a_n$ then $(x,x,\widetilde{e}_0)=0$ because $\ache_x$ is associative.

 For $x+r\widetilde{e}_0$ with $x$ doubly pure we have
\begin{eqnarray*}
(x+r\widetilde{e}_0,x+r\widetilde{e}_0,\widetilde{e}_0)
&=&
(x,x,\widetilde{e}_0)+r(x,\widetilde{e}_0,\widetilde{e}_0)+r(\widetilde{e}_0,
x,\widetilde{e}_0)+r^2(\widetilde{e}_0,\widetilde{e}_0,\widetilde{e}_0)\\
&=&0
\end{eqnarray*}
because
 $\widetilde{e}_0$ is alternative.Therefore  $\widetilde{e}_0$ is strongly alternative. 

We will show that the only strongly alternative elements are of the form $re_0+s\widetilde{e}_0$ for $r$ and $s$ in $\erre$.
\vglue.5cm
\noindent
{\bf Lemma 4.1.} Let $a$ be a non-zero element in $\a_n$ for $n\geq 3$.

 If 
$(x,a,y)=0$ for all $x$ and $y$ in $\a_n$ then $a=re_0$ for $r\in\erre$.
\vglue.5cm
\noindent
{\bf Proof.} First of all we notice that $a$ has to be 
an
alternative element because if $x=a$ then $(a,a,y)=0$ for all $y$ in $\a_n$.

Now $(x,a,y)=0$ for all $x$ and $y$ if and only if the pure part of $a$ also
has this property, i.e. write $a=b+re_0$ for $b$ pure element and $r\in\erre$ 
then
 $(x,a,y)=(x,b,y)+(x,re_0,y)$ but $(x,re_0,y)=0$ for all $x$ and $y$ in $\a_n$ and $r$ in $\erre$. Since $b$ is also alternative its doubly pure part is  alternative so if $b=c+s\widetilde{e}_0$ with $c$ doubly pure and $s$ in $\erre$ then $c$ is alternative. 
Setting $x=\widetilde{e}_0$ we have that
$$0=(\widetilde{e}_0,b,y)=(\widetilde{e}_0,c+s\widetilde{e}_0,y)=(\widetilde{e}_0,c,y)+s(\widetilde{e}_0,\widetilde{e}_0,y)$$%
$$=-2\widetilde{c}y+0$$
for all $y$ in $\a_n$. 
But $c$ is alternative so $c=0$. Therefore $b=s\widetilde{e}_0$ with $s$ in $\erre$. By hypothesis $(x,b,y)=s(x,\widetilde{e}_0,y)=0$ for all $x$ and $y$, so $s=0$ and $b=0$ in $\a_n$. So if $(x,a,y)=0$ for all $x$ and $y$ in $\a_n$ then $a$ has imaginary part equal to zero.

\hfill Q.E.D.
\vskip.3cm
\noindent
{\bf Theorem 4.2.} If $\alpha\in\a_{n+1}\quad n\geq 3$ is a pure strongly alternative element then $\alpha=\lambda\widetilde{e}_0$ for some $\lambda\in
\erre$.
\vskip.3cm
\noindent
{\bf Proof.} First suppose that $\alpha$ is doubly pure, and strongly alternative in $\a_{n+1}$.
We want to show that $\alpha=0$.

Using theorem 3.3, since $\alpha$ is alternative then $\alpha=(ra,ta)\in\a_n\times\a_n$ with $a$ an alternative pure element in $\a_n$ and $r$ and $t$ fixed real numbers, not both equal to zero (otherwise we are done with $\alpha=0$). Now for all 
$\gamma=(x,y)\in _o\!\!\a_n\times_o\!\!\a_n$ we have that
$$L^2_\gamma(\alpha)=({\cal A}(ra)-S(ta),{\cal A}(ta)+S(ra))$$
where ${\cal A}=L^2_x+R^2_y$ and $S=(x,-,y)$ (see $\S 3$). Since $\alpha$ is strongly alternative then
$$L^2_\gamma(\alpha)=\gamma^2\alpha\quad{\rm for  all}\quad 
\gamma\in\widetilde{\a}_{n+1}$$
thus we have the following system of equations
\begin{eqnarray*}
r{\cal A}(a)-tS(a)&=&(x^2+y^2)ra\\
t{\cal A}(a)+rS(a)&=&(x^2+y^2)ta.
\end{eqnarray*}
If $r=0$ or $t=0$ then $S(a)=0$. If $r\neq 0$ and $t\neq 0$ then also $S(a)=0$ because multiplying the first equation by `$t$', the second by `$r$' and substracting we have that
$$(r^2+t^2)S(a)=0.$$
Therefore $S(a)=0$. But $a$ is pure element so by Lemma 4.1 $a=0$ in $\a_n$ and $\alpha=0$ in $\a_{n+1}$.

Suppose now that $\alpha=\beta+\lambda\widetilde{e}_0$ with $\beta$ double 
pure  in $\a_{n+1}$ and $\lambda$ in $\erre$. Since $(\alpha,x,x)=0$ for all $x$ in $\a_{n+1}$ then
$0=(\beta+\lambda\widetilde{e}_0,x,x)=(\beta,x,x)+\lambda(\widetilde{e}_0,x,x)$ but $(\widetilde{e}_0,x,x)=0$ for all $x$ (see example above). 
Therefore $0=(\beta,x,x)$ for all $x$ in $\a_{n+1}$ and $\beta=0$ so $\alpha=\lambda\widetilde{e}_0$ for $\lambda$ in $\erre$.

\hfill Q.E.D.
\vskip.3cm
\noindent
{\bf Corollary 4.3.} If $\varphi\in Aut(\a_n)$ is an algebra authomorphism of $\a_n$ for $n\geq 4$ then $\varphi(\widetilde{e}_0)=\pm\widetilde{e}_0$ and 
$\varphi(\widetilde{a})=\pm\widetilde{\varphi(a)}$ for all $a\in\a_n$.
\vglue.5cm
\noindent
{\bf Proof.} If $\varphi\in Aut(\a_n)$ then $\varphi(e_0)=e_0$ and $\varphi$ preserves real parts.

Now $\varphi$ preserves the properties of being alternative and strongly alternative so $\varphi(\widetilde{e}_0)=\lambda\widetilde{e}_0$ for some $\lambda$ in $\erre$ and $\varphi\in Aut(\a_n)$ but
$$1=||\varphi(\widetilde{e}_0)||^2=||\lambda\widetilde{e}_0||^2=\lambda^2||\widetilde{e}_0||=\lambda^2\quad{\rm so}\quad \lambda=\pm1.$$
\hfill Q.E.D.
\vskip.3cm
\noindent
{\bf Remark.} From the last corollary we can deduce the group of authomorphism of $\a_{n+1}$ in terms of the group of authomorphism of $\a_n$ for $n\geq 3$.

For $\varphi\in Aut(\a_{n+1}), \varphi$ is completly determined by the action on ``the even part of $\a_{n+1}$'' i.e. $\{(a,0)\in\a_n\times\a_n\}$ which is isomorphic
 to $\a_n$ and the sign $\varphi(\widetilde{e}_0)=\pm\widetilde{e}_0$.

Recall that $(a,b)=(a,0)+(b,0)(0,e_0)$ in $\a_n\times\a_n$.

 Eakin-Sathaye already calculate this group [1]
$$Aut(\a_{n+1})=Aut(\a_n)\times\sum_3\qquad n\geq 3$$
where $\sum_3$ is the symmetric group with 6 elements given by $\langle\tau,\mu:\tau^2=\mu^3=1$ and $\mu\tau=\tau\mu^2\rangle , \tau:\a_{n+1}\rightarrow\a_{n+1}$ is $\tau(x,y)=(x,-y)$ then $\tau(\widetilde{e}_0)=-\widetilde{e}_0$ and $\mu(x,0)=(x,0)$ and $\mu(0,x)=(x,0)\alpha$ where $\alpha=-{1\over 2}e_0-{1\over 2}\widetilde{e}_0$.
\vglue1cm
\noindent
{\bf \S 5 Local alternativity.}
\vglue.5cm
\noindent
{\bf Definition.} $a$ and $b$ are in $\a_n$ with $n\geq 4$. We say that 
$a$ {\it is normed with} $b$ if
$$||ab||=||a||||b||.$$
We say that $a$ {\it alternate with} $b$ if
$$(a,a,b)=0.$$
Way say that $a$ alternate strongly with $b$ if $(a,a,b)=0$ and $(a,b,b)=0$.

Now $a$ is {\it a normed element} if $||ax||=||a||||x||$ for all $x$ in $\a_n$.

Notice that since `$e_0$' is normed, alternative and strongly alternative element we have that $a$ is normed with $b$, $a$ alternate with $b$ and $a$ alternate strongly with $b$ if the pure part of $a$ is normed, alternate and strongly 
alternate with b.

Now if $a$ is a pure element in $\a_n$ with $n\geq 4$ then $a$ alternative 
with $b$ implies that $a$ is normed with $b$ for $b$ in $\a_n$.

To see this we argue as follows:
\begin{eqnarray*}
||ab||^2&=&\langle ab,ab\rangle=-\langle b,a(ab)\rangle=-\langle b,a^2b\rangle=-a^2\langle b,b\rangle\\
&=&||a||^2||b||^2\quad{\rm for}\quad b\quad{\rm in}\quad \a_n.
\end{eqnarray*}
The converse is not true.
\vglue.5cm
\noindent
{\bf Example.} $a=e_1+e_{10}$ and $b=e_{15}$ in $\a_4$.
\begin{eqnarray*}
(a,a,b)&=&(e_1,e_1,e_{15})+(e_{10},e_{10},e_{15})+(e_1,e_{10},e_{15})+(e_{10},e_1,e_{15})\\
&=&0+0+[(e_1e_{10})+(e_{10}e_1)]e_{15}-(e_1(e_{10}e_{15})+e_{10}(e_1e_{15}))\\
&=&0+0+0-(e_1e_5-e_{10}e_{14})\\
&=&-(-e_4-e_4)\\
&=&2e_4.
\end{eqnarray*}
On the other hand
\begin{eqnarray*}
||ab||^2&=&||(e_1+e_{10})e_{15}||^2=|e_1e_{15}+e_{10}e_{15}||^2=||
-e_{14}+e_5||^2\\
&=&||e_5||^2+||e_{14}||^2-2\langle e_5,e_{14}\rangle\\
&=&2\\
||a||^2||b||^2&=&||e_1+e_{10}||^2||e_{15}||^2=||e_1||^2+||e_{10}
||^2+2\langle e_1,e_{10}\rangle\\
&=&2
\end{eqnarray*}
so $||ab||=||a||||b||$.

Notice that if $||ab||=||a||||b||$ with $a$ pure  then
$$-\langle b,a^2b\rangle=-a^2\langle b,b\rangle=||a||^2||b||^2=||ab||^2=\langle ab,ab\rangle=-\langle b,a(ab))\rangle .$$
Therefore $\langle b,(a,a,b)\rangle=0$ and $(a,a,b)\perp b$ 
which is equivalent to $(a,b,b)\perp a$. Also if $a$ alternate strongly with 
$b$ then $a$ alternate with $b$ (trivially) but the converse is also 
not true. (See example in $\S 4$). Therefore

\hglue2cm ``$a$ alternate strongly with $b$''

\hglue4cm $\Downarrow$

\hglue2cm ``$a$ alternate  with $b$''

\hglue4cm $\Downarrow$

\hglue2cm ``$a$ is normed with $b$''

\noindent
and the converse of both implications are not true. Notice that since
$\a_3$ is an alternative normed algebra the three concepts are equivalent in
$\a_3$ 

Since $\widetilde{e}_0$ in $\a_n ,\quad n\geq 4$, is a
strongly alternative element, we have that $a$ alternate strongly with $b$ if the doubly pure part of $a$ strongly alternate with $b$ and viceversa. 
So we are interested in these relations between couples of doubly pure elements in $\a_n$ for $n\geq 4$. 

Now for a (non-zero) doubly pure element we have that if $b\in\ache_a$ then $a$ alternate strongly with $b$ because $\ache_a$ is associative.
\vglue.5cm
\noindent
{\bf Theorem 5.1.} 

Let $a$ and $b$ be (non-zero) doubly pure elements 
in $\a_n$ with $n\geq 4$.

If $b\in\ache^\perp_a$ and $a$ alternate strongly with $b$ of $\a_n$.

Then
i) The vector subspace of $\a_n$ generated by
$$\{e_0,a,b,ab\}.$$
is multiplicatively closed and isomorphic to $\a_2=\ache$.

ii) The vector subspace of $\a_n$ generated by
$$\{e_0,a,b,ab,\widetilde{a}b,-\widetilde{b},\widetilde{a},\widetilde{e}_0\}$$
is multiplicatively closed and isomorphic to $\a_3=\o$.
\vglue.5cm
\noindent
{\bf Proof.} Without loss of  generality we assume that
$||a||=||b||=1$ (otherwise we take ${a\over||a||}$ and ${b\over ||b||})$.

Now construct a multiplication table
$$
\begin{tabular}{l|llll}
&$e_0$&$a$&b&$ab$\\ \hline
$e_0$&$e_0$&$a$&$b$&$ab$\\
$a$&$a$&$-e_0$&$ab$&$-b$\\
$b$&$b$&$-ab$&$-e_0$&$a$\\
$ab$&$ab$&$b$&$-a$&$-e_0$
\end{tabular}
$$
Since $b\perp a, ab=-ba$ and $(a,a,b)=(a,b,b)=0$ we have that 
$a(ab)=a^2b=-b; b(ab)=-b(ba)=-b^2a=a; (ab)a=-(ba)a=-ba^2=b; (ab)b=ab^2=-a$ and
$(ab)^2=-||ab||^2=-\langle ab,ab\rangle=\langle b,a(ab)\rangle=\langle b,a^2b\rangle=a^2||b||^2=-e_0$.

This multiplication table is the one of $\a_2$ by the identification 
$e_0\leftrightarrow e_0; e_1\leftrightarrow a; e_2\leftrightarrow b; 
 e_3\leftrightarrow ab$. So we are done with i).

ii) is a routine calculation using i) and Lemma 1.1 and Corollary 1.2.

\hfill Q.E.D.
\vskip.3cm
\noindent
{\bf Remarks:} Given a non-zero element $a$ in $\a_n\quad n\geq 4$ there exists $b$ such that: $a\perp b, ||a||=||b||$ and $a$ alternate strongly with $b$.

To see this take $a$ of normed one (otherwise we take ${a\over||a||}$) and write $a=rc+s\widetilde{e}_0$ where $c$ is the doubly pure part of $a$ and $r$ and $s$ are in $\erre$ with $r^2+s^2=1$.

Define $b=sc-r\widetilde{e}_0$. It is a routine $tb$ see that $(a,a,b)=(a,b,b)=$ and that $||a||=||b||$ and $a\perp b$.
\vskip.3cm
\noindent
{\bf Definition.} Let $B$ be a subset of $\a_n\quad n\geq 4\quad B$ is a normed subset of $\a_n$ if $||xy||=||x||||y||$ for all $x$ and $y$ in $B$. Thus $a$ is normed with $b$ means that $\{a,b\}$ is a normed subset of $\a_n$.
\vskip.3cm
\noindent
{\bf Theorem 5.2} Let $a$ and $b$ be non-zero pure elements in $\a_n\quad n\geq 4$.

1) $a$ alternate with $b$ if and only if $\{a,b\}$ and $\{a,ab\}$ are normed subsets.

2) $a$ alternate strongly with $b$ if and only if $\{a,b,ab\}$ is a normed subset.
\vskip.3cm
\noindent
{\bf Proof.} 1) We know that if $(a,a,x)=0$ then $||ax||=||a||||x||$ for all $x$ in $\a_n$.

But $a(a,a,b)=a(a^2b-a(ab))=a^2ab-a(a(ab))=(a,a,ab)$ so $(a,a,b)=0$ if and only if $(a,a,ab)=0$ then $\{a,b\}$ and $\{a,ab\}$ are normed subsets when $a$ alternate with $b$.

Conversely suppose that $||ab||=||a||||b||$ and that $||a(ab)||=||a||||ab||=||a||^2||b||$. Then
\begin{eqnarray*}
||(a,a,b)||^2&=&\langle(a,a,b),(a,a,b)\rangle\\
&=&\langle a^2b-a(ab), a^2b-a(ab)\rangle\\
&=&a^4\langle b,b\rangle+\langle a(ab),a(ab)\rangle-2\langle a^2b, a(ab)\rangle\\
&=&||a||^4||b||^2+||a(ab)||^2+2||a||^2\langle +b, a(ab)\rangle\\
&=&||a||^4||b||^2+||a||^4||b||^2-2||a||^2\langle ab,ab\rangle\\
&=&2||a||^4||b||^2-2||a||^2||ab||^2\\
&=&0
\end{eqnarray*}
Therefore $(a,a,b)=0$ and we are done with 1). To prove 2) we apply 1) to $\{a,b\}, \{a,ab\}$ and $\{b,ba\}$ but $||ba||=||ab||$ so $(a,a,b)=(a,b,b)=0$ if and only if $\{a,b,ab\}$ is a normed set.

\hfill Q.E.D.

\newpage
\noindent
{\bf References}
\vglue.5cm
\begin{enumerate}
\item[{\rm [1]}]P. Eakin - A. Sathaye. On the automorphisms of the Cayley-Dickson algebras and its derivations. Journal of algebra 129, 263-278 (1990).
\item[{\rm [2]}] S.H. Khalil-P. Yiu. The Cayley-Dickson algebras. Bolet\'{\i}n de la Sociedad de Lodz, Vol. XLVIII 117-169, 1997.
\item[{\rm [3]}]G. Moreno. The zero divisors of the Cayley-Dickson algebras over the real numbers. Boletin de la Sociedad Matem\'atica Mexicana (3), Vol. 4, 13-27, 1998.
\item[{\rm [4]}]R.D. Schafer. On the algebras formed by the Cayley-Dickson process. American Jorunal of Math. 76, 1954, 435-446.
\end{enumerate}

\end{document}